\documentclass[11pt,reqno]{amsart}
\usepackage{amsmath,amssymb,bbm,times,url}

\newtheorem{theorem}{Theorem}
\newtheorem{lemma}[theorem]{Lemma}
\newtheorem{proposition}[theorem]{Proposition}
\newtheorem{corollary}[theorem]{Corollary}

\newcommand{\e}{\varepsilon}
\newcommand{\charfn}{\mathbbm{1}}
\newcommand{\cube}{{\mathcal C}}
\newcommand{\C}{{\mathbb C}}

\newcommand{\Imag}{\operatorname{Im}}

\newcommand{\N}{{\mathbb N}}

\newcommand{\R}{{\mathbb R}}
\newcommand{\Rd}{{\R^d}}
\newcommand{\Real}{\operatorname{Re}}
\newcommand{\Rem}{\operatorname{Rem}}

\newcommand{\spt}{\operatorname{spt}}
\newcommand{\Z}{{\mathbb Z}}
\newcommand{\Zd}{{\Z^d}}

\begin{document}

\title[]{Affine synthesis onto $L^p$ when $0 < p \leq 1$}

\author[]{R. S. Laugesen}
\address{Department of Mathematics, University of Illinois, Urbana,
IL 61801, U.S.A.} \email{Laugesen\@@uiuc.edu}
\date{\today}

\thanks{Laugesen's travel was supported by N.S.F. Award DMS--0140481.}

\begin{abstract}
The affine synthesis operator $Sc=\sum_{j>0} \sum_{k \in \Zd}
c_{j,k} \psi_{j,k}$ is shown to map the coefficient space
$\ell^p(\Z_+ \times \Zd)$ surjectively onto $L^p(\Rd)$, for $p \in
(0,1]$. Here $\psi_{j,k}(x)=|\det a_j|^{1/p} \psi(a_j x - k)$ for
dilation matrices $a_j$ that expand, and the synthesizer $\psi \in
L^p(\Rd)$ need satisfy only mild restrictions, for example $\psi \in
L^1(\Rd)$ with nonzero integral or else with periodization that is
real-valued, nontrivial and bounded below.

An affine atomic decomposition of $L^p$ follows immediately:
\[
\| f \|_p \approx \inf \left\{ ( \sum_{j>0} \sum_{k \in \Zd}
|c_{j,k}|^p )^{1/p} : f=\sum_{j>0} \sum_{k \in \Zd} c_{j,k}
\psi_{j,k} \right\} .
\]

Tools include an analysis operator that is nonlinear on $L^p$.
\end{abstract}

\maketitle

\section{\bf Introduction}
\label{introduction}

Many normed function spaces can be generated by discrete translates
and dilates of just a single function. For example, Sobolev spaces
can be decomposed by spline approximation or wavelet expansion. But
in metric vector spaces that are not normed, the theory of such
affine systems is much less developed. This paper develops the
affine theory of $L^p=L^p(\Rd), 0 < p \leq 1$.

Given a synthesizer $\psi \in L^p$, the \emph{affine synthesis}
operator is
\[
c = \{ c_{j,k} \} \mapsto \sum_{j > 0} \sum_{k \in \Zd} c_{j,k}
\psi_{j,k} = Sc
\]
where
\[
\psi_{j,k}(x)=|\det a_j|^{1/p} \psi(a_j x - bk) , \qquad x \in \Rd
.
\]
The \emph{dimension} $d \in \N$ and the exponent $p \in (0,1]$ are
fixed. The \emph{dilation matrices} $a_j$ are invertible $d \times
d$ real matrices that are \emph{expanding}, in the sense that their
inverses contract to zero:
\[
\| a_j^{-1} \| \to 0 \qquad \text{as $j \to \infty$}
\]
where $\| \cdot \|$ denotes the operator norm of the matrix acting
from $\Rd$ to $\Rd$. For example, one could take $a_j=2^j I$. The
\emph{translation matrix} $b$ is an invertible $d \times d$ real
matrix, for example the identity. Note this paper only uses $j>0$,
meaning the affine systems only use small scales.

Our first goal is to find the right domain for the synthesis
operator, that is, to find a sequence space that $S$ maps
continuously into $L^p$. Proposition~\ref{continuitysynth} shows $S$
is continuous from $\ell^p(\Z_+ \times \Zd)$ to $L^p$, where $\Z_+ =
\{ j \in \Z : j>0 \}$.

Our second goal is surjectivity: we want $S$ to map $\ell^p(\Z_+
\times \Zd)$ \emph{onto} $L^p$, so that every function in $L^p$ can
be written as an infinite linear combination of the $\psi_{j,k}$.
Theorem~\ref{surjectivity} proves this surjectivity, by building on
H\"{o}lder continuity of the nonlinear analysis operator as
introduced in Theorem~\ref{continuityanal}, and using an explicit
$L^p$-approximation result in Theorem~\ref{sample}. The underlying
idea, roughly, is to quasi-interpolate via nonlinear analysis and
then linear synthesis, at a very small scale, and then to apply the
open mapping theorem.

To illustrate our result, observe that Theorem~\ref{surjectivity}
and the sufficient condition in Proposition~\ref{bitten} combine to
yield that if $\psi \in L^p \cap L^1, p \in (0,1)$, and either
$\int_\Rd \psi \, dx \neq 0$ or else $0 \not \equiv \sum_{k \in \Zd}
\psi(x-bk)$ is real-valued and bounded below, then $S$ maps
surjectively onto $L^p$. Thus surjectivity holds for a large class
of synthesizers $\psi$. Indeed, Theorem~\ref{properties} takes a
global perspective and shows that surjectivity of the synthesis
operator holds generically with respect to the choice of synthesizer
$\psi \in L^p$.

Interestingly, Strang--Fix conditions are not required in this
paper: the integer translates of the synthesizer need not form a
partition of unity. But if these translates do sum up to $1$, then
the bounds in the surjectivity result Theorem~\ref{surjectivity} get
better (because one can take $\sigma=0$ and $\lambda=1$ there).

Surjectivity of $S$ onto $L^p$ immediately implies an affine atomic
decomposition (or metric equivalence) in Corollary~\ref{atomic}, of
the form
\[
\| f \|_p \approx \inf \{ \| c \|_{\ell^p(\Z_+ \times \Zd)} : f=Sc
\} .
\]
When $p=1$ this was found earlier by Bruna \cite[Theorem~4]{B06}.
Corollary~\ref{atomicdomain} localizes the atomic decomposition to
$L^p(\Omega)$, for domains $\Omega \subset \Rd$.

Theorem~\ref{inject} restricts to a single dilation scale $j$, and
states an atomic decomposition that does not need an ``$\inf$''. In
other words, it proves when $\psi$ has compact support that the
$\psi_{j,k}, k \in \Zd$, form a $p$-Riesz basis for their closed
linear span in $L^p$, or that the synthesis operator at scale $j$ is
bounded, injective and has closed range in $L^p$. This result
slightly extends some of Jia's work \cite{J98} on $L^p$-stability of
shift invariant subspaces.

The results of this paper for $0 < p \leq 1$ are contrasted with
prior work on $p \geq 1$ in Section~\ref{remarksLHS}. Open problems
are raised in Section~\ref{openproblems}, including Meyer's Mexican
hat spanning problem for $L^p$ when $1<p<\infty$, which this paper
resolves for $0<p<1$.

\subsubsection*{Discussion.} This paper shows that
arbitrary $L^p$ functions can be decomposed into linear combinations
of discrete translates and dilates of the synthesizer $\psi$,
without requiring any particularly special properties of $\psi$.
This structural information about $L^p$ has intrinsic mathematical
interest, and might conceivably be useful in applications for which
a particular shape of $\psi$ is naturally preferred.

The central contribution of the paper is its constructive method of
$L^p$-controlled approximation via nonlinear analysis and linear
synthesis (Theorem~\ref{sample}), which implies surjectivity of the
synthesis operator (Theorem~\ref{surjectivity}). The closest prior
result for $p \in (0,1)$ is due to Filippov and Oswald
\cite{FO95,F98}, who proved for isotropic dilation matrices that
every $L^p$ function can be written as $Sc$ for some sequence $c$,
but unfortunately with no information on the size of $c$ or to what
space $c$ might belong.

Incidentally, DeVore \emph{et al.}\ \cite{DJP92} have proved that
linear combinations of the $\psi_{j,k}$ provide good
$L^p$-approximations to functions in Besov spaces, and an abstract
framework for that was developed in \cite{CT97}.

\subsubsection*{Notation.} $L^p=L^p(\Rd)$ denotes the
class of complex valued functions with $\| f \|_p = (\int_\Rd |f|^p
\, dx)^{1/p} < \infty$. It is a complete metric space with distance
function
\[
d_p(f,\widetilde{f}) = \| f - \widetilde{f} \|_p^p .
\]
Write $f \equiv 0$ to mean $f=0$ a.e., that is $\| f \|_p = 0$.

A multi-scale, discrete analogue of $L^p$ is the space $\ell^p(\Z_+
\times \Zd)$ consisting of doubly-indexed sequences $c=\{ c_{j,k}
\}_{j>0,k \in \Zd}$ of complex numbers satisfying
\[
\| c \|_{\ell^p(\Z_+ \times \Zd)} = ( \sum_{j>0} \sum_{k \in \Zd}
|c_{j,k}|^p )^{1/p} < \infty .
\]
Clearly $\ell^p(\Z_+ \times \Zd)$ is a complete metric space with
distance function
\[
d_{\ell^p}(c,\widetilde{c}) = \| c - \widetilde{c} \|_{\ell^p(\Z_+
\times \Zd)}^p .
\]

\subsubsection*{Useful fact.} The triangle inequality
for the $\ell^p$-metric on the complex numbers says
\[
|\sum_m z_m|^p \leq \sum_m |z_m|^p , \qquad z_m \in \C , \quad p \in
(0,1] .
\]
In other words, $p$-th powers can be taken inside sums.

\section{\bf Results} \label{results}

Our first four results show that synthesis maps continuously into
$L^p$, that nonlinear analysis is continuous on $L^p$, that
synthesis and analysis can partially reconstruct every $L^p$
function, and hence that synthesis maps surjectively onto $L^p$.
Next we deduce an affine atomic decomposition of $L^p$, and a metric
equivalence via the analysis operator. Then we prove synthesis at
each fixed dilation scale gives a $p$-Riesz basis. Our last result
considers the class of all synthesizers for which the synthesis
operator is surjective, and proves the class is dense, open and
connected in $L^p$.

Recall from the introduction that the \emph{synthesis operator} is
\begin{equation} \label{Sdef}
Sc = \sum_{j>0} \sum_{k \in \Zd} c_{j,k} \psi_{j,k} .
\end{equation}
Synthesis is Lipschitz continuous:
\begin{proposition}[Synthesis into $L^p$] \label{continuitysynth}
Assume $\psi \in L^p, p \in (0,1]$.

Then $S : \ell^p(\Z_+ \times \Zd) \to L^p$ is continuous and linear.
More precisely, if $c \in \ell^p(\Z_+ \times \Zd)$ then the series
\eqref{Sdef} for $Sc$ converges pointwise absolutely a.e.\ to a
function in $L^p$ (and hence $Sc$ converges unconditionally in
$L^p$), and
\begin{equation} \label{analogue1}
\| Sc \|_p \leq \| \psi \|_p \| c \|_{\ell^p(\Z_+ \times \Zd)} .
\end{equation}
\end{proposition}
The proof is in Section~\ref{continuitysynth_proof}.

\vspace{6pt} Next we develop our analysis operator. The traditional
linear definition of analysis is clearly invalid on $L^p$ for $p<1$,
for if we tried to integrate $f \in L^p$ against an analyzing
function then the integral might not even exist, since $f$ need not
be locally integrable. We compensate for this lack of local
integrability by applying a nonlinear radial stretch to $f$ before
analyzing it, and undoing the stretch afterwards.

The radial stretch function
\[
\Theta(z) = |z|^{p-1} z, \qquad z \in \C ,
\]
is a homeomorphism of the complex plane and satisfies
$|\Theta(z)|=|z|^p$. It acts on complex-valued functions $f$ by
\[
(\Theta f)(x) = \Theta(f(x)) .
\]
Notice that if $f \in L^p$ then $\Theta f \in L^1$.

Define the \emph{analysis operator at scale $j$} by
\[
T_j f = \{ |\det b| \Theta^{-1} \langle \Theta f , \phi(a_j \,
\cdot - bk) \rangle \}_{k \in \Zd} .
\]
Roughly, $T_j$ maps a function $f$ to its sequence of sampled
$\phi$-averages at scale $j$, except that radial stretching is
applied to $f$ before the sampling, and then is undone at the end.
The $|\det b|$ factor is for later convenience. We emphasize that
\begin{quote}
the analysis operator $T_j$ is \emph{nonlinear}, and depends
implicitly on the exponent $p$.
\end{quote}
(The synthesis operator also depends on $p$, through the
normalization of $\psi_{j,k}$.)

The next theorem shows the analysis operator is continuous. The
hypotheses involve the \emph{periodization} operator, defined on a
function $g$ by
\[
Pg(x) = |\det b| \sum_{k \in \Zd} g(x-bk) \qquad \text{for $x \in
\Rd$.}
\]
Clearly $Pg$ is periodic with respect to the lattice $b\Zd$,
provided it is well-defined a.e.
\begin{theorem}[Analysis on $L^p$] \label{continuityanal}
Assume $p \in (0,1]$ and take $\phi \in L^\infty$ with $P|\phi|
\in L^\infty$.

Then for each $j$, the analysis operator $T_j : L^p \to \ell^p(\Zd)$
is locally H\"{o}lder continuous, with
\[
d_{\ell^p}(T_j f , T_j g) \leq C [d_p(f,0) + d_p(0,g)]^{1-p}
d_p(f,g)^p , \qquad f, g \in L^p .
\]
Here $C$ depends on the translation matrix $b$, the exponent $p$
and the analyzer $\phi$, but not on the dilation scale $j$.
\end{theorem}
Section~\ref{continuityanal_proof} has the proof. Recall the
distance function on $\ell^p(\Zd)$ is $d_{\ell^p}(s,t)=\| s - t
\|_{\ell^p(\Zd)}^p$ when $p \in (0,1]$. The hypothesis that the
absolute value $|\phi|$ of the analyzer have bounded periodization
is easily satisfied, say if $\phi$ is bounded with compact support
or with rapid decay.

The nonlinear analysis operator $T_j$ is also locally H\"{o}lder
continuous when $p>1$, with
\[
d_{\ell^p}(T_j f , T_j g) \leq C [d_p(f,0) + d_p(0,g)]^{1-1/p}
d_p(f,g)^{1/p} , \qquad f, g \in L^p ,
\]
where for $p>1$ the distance functions are defined by
$d_{\ell^p}(s,t)=\| s - t \|_{\ell^p(\Zd)}$ and $d_p(f,g)=\| f - g
\|_p$. We omit the proof, since the rest of this paper concerns $p
\in (0,1]$.

\vspace{6pt} Now we start to develop approximation results. Write
\[
S_j s = \sum_{k \in \Zd} s_k \psi_{j,k}
\]
for the \emph{synthesis operator at scale $j$}, acting on sequences
$s=\{ s_k \}_{k \in \Zd}$. Notice $S_j$ is continuous from
$\ell^p(\Zd)$ to $L^p$, with
\[
\| S_j s \|_p \leq \| \psi \|_p \| s \|_{\ell^p(\Zd)}
\]
for $p \in (0,1]$ by Proposition~\ref{continuitysynth}.

The following approximation result will be used later in proving
surjectivity of the synthesis operator. It uses
\[
\cube=[0,1)^d
\]
to denote the unit cube in $\Rd$.
\begin{theorem}[Affine quasi-interpolation] \label{sample}
Assume $\psi \in L^p, p \in (0,1]$, and suppose $\phi \in L^\infty$
with $P|\phi| \in L^\infty$ and $\int_\Rd \phi \, dx = 1$.

Then
\[
\lim_{j \to \infty} \| S_j T_j f - f \|_p = \| P\psi(b \, \cdot) - 1
\|_{L^p(\cube)} \| f \|_p , \qquad f \in L^p .
\]
In particular, if $\psi$ has constant periodization $P \psi = 1$
a.e., then $S_j T_j f \to f$ in $L^p$ as $j \to \infty$.
\end{theorem}
Section~\ref{sample_proof} has the proof. Note the periodization
series $P\psi(bx)$ appearing in the theorem is well defined whenever
$\psi \in L^p$, because it converges absolutely a.e.\ and belongs to
$L^p(\cube)$:
\begin{align*}
\int_\cube |P\psi(bx)|^p \, dx & \leq |\det b|^p \int_\cube \sum_{k
\in \Zd} |\psi(b(x-k))|^p \, dx \qquad \text{using $p \in (0,1]$} \\
& = |\det b|^p \| \psi(b \, \cdot) \|_p^p < \infty .
\end{align*}

The constant periodization condition $P \psi = 1$ says that the
collection $\{ |\det b| \psi(x - bk) : k \in \Zd \}$ of translates
of $\psi$ is a partition of unity. Examples of such $\psi$ (when
$b=I$) include the indicator function $\charfn_\cube$ and
$B$-splines obtained by convolution with this indicator function.

When $p \in [1,\infty)$ and $\psi$ has constant periodization, the
result that $S_j T_j f \to f$ (with $T_j$ being a linear analysis
operator) has a long history in Strang--Fix approximation theory,
summarized in \cite[\S3]{bl2}.

%The compact support assumption on the analyzer $\phi$, in
%Theorem~\ref{sample}, can be weakened to $\phi \in L^\infty$ being
%bounded by a radially decreasing $L^1$ function. Indeed one can
%prove the map $\phi \mapsto T_j f$ is locally $p$-H\"{o}lder
%continuous with respect to the norm $\| P|\phi| \|_\infty$
%(uniformly in $j$ and for fixed $f$), which enables a reduction to
%the compactly supported case of the theorem. But we will not need
%this extension.

\vspace{6pt} Next we show every $f \in L^p$ can be written as $Sc$
for some sequence $c \in \ell^p(\Z_+ \times \Zd)$, so that the
synthesis operator is surjective.
\begin{theorem}[Synthesis onto $L^p$] \label{surjectivity}
Assume $\psi \in L^p, p \in (0,1]$, and suppose that
\begin{equation} \label{bite}
\sigma := \| \lambda P\psi(b \, \cdot) - 1 \|_{L^p(\cube)}^p < 1
\end{equation}
for some $\lambda \in \C$.

Then $S : \ell^p(\Z_+ \times \Zd) \to L^p$ is open, and surjective.
Indeed, if $f \in L^p$ and $\sigma^\prime \in (\sigma,1)$ then a
sequence $c \in \ell^p(\Z_+ \times \Zd)$ exists such that $Sc=f$ and
\[
\| c \|_{\ell^p(\Z_+\times \Zd)} \leq (1-\sigma^\prime)^{-1/p}
|\lambda| |\det b|^{1-1/p} \| f \|_p .
\]
\end{theorem}
The proof is in Section~\ref{surjectivity_proof}. We discuss the
hypothesis \eqref{bite} below.

We do not know any prior general work on surjectivity of the
synthesis operator when $p \in (0,1)$. The closest seems to be
Filippov and Oswald's construction in \cite{FO95,F98} of
``representation systems'', by which every $f \in L^p$ can be
written as a convergent series $Sc=f$, provided the dilation
matrices are real multiples of the identity. This looks like
surjectivity, but the drawback is that their result yields no
control over the size of coefficients in the sequence $c$, and thus
it is unclear what the domain of the synthesis operator actually is.
The achievement of Theorem~\ref{surjectivity} is to identify the
sequence space $\ell^p(\Z_+ \times \Zd)$ as a domain from which $S$
maps onto $L^p$.

Filippov and Oswald's method involves iterative approximation of
simple functions, and thus is less concrete than our approach. Their
paper uses the open mapping theorem implicitly, whereas we use it
explicitly. Interestingly, their key lemma is crucial to our proof
too (Proposition~\ref{bitten} below), although we employ it
differently.

\subsubsection*{Non-injectivity of $S$.} The synthesis
operator is certainly not injective, and has a very large kernel.
For example, we could discard the dilation $a_1$ (discarding all
terms with $j=1$ in the sum defining $Sc$) and still show $S$ maps
onto $L^p$, by applying Theorem~\ref{surjectivity} with the
remaining dilations $\{ a_2,a_3,\ldots \}$.

\vspace{6pt} For Theorem~\ref{surjectivity} to be interesting, we
need to exhibit examples of synthesizers $\psi$ satisfying
hypothesis \eqref{bite}.
\begin{quotation}
\noindent {\small \textsc{Examples.} If $\psi$ is supported in the
cube $\cube$ with $\psi \geq 1$ there, and with $\| \psi - 1
\|_{L^p(\cube)} < 1$ and $b=I$, then $\psi$ equals its own
periodization on $\cube$, and so \eqref{bite} holds with $\lambda =
1$. To be specific, in one dimension one could take
$\psi(x)=(1+Ax^{-\beta}) \charfn_{[0,1)}$ for $\beta<1/p$ and
suitably small $A>0$. This example shows $\psi$ can have a typical
$L^p$-singularity at the origin.

Another singular example in one dimension is $\psi(x) = x^{-\beta}
\charfn_{[0,1)}$ for $\beta \in (0,2/(p+1))$, which we prove in
Section~\ref{tachev_proof} satisfies \eqref{bite} when $b=1$.}
\end{quotation}
On the other hand, some functions $\psi$ do not satisfy
\eqref{bite}:
\begin{quotation} \noindent {\small \textsc{Counterexample.} Tachev
\cite{T95} proved that \eqref{bite} with $b=1$ fails in one
dimension for $\psi(x) = x^{-\beta} \charfn_{[0,1)}$ when $\beta \in
[2/(p+1),1/p)$ and $p \in (0,1)$. We show this in
Section~\ref{tachev_proof}.}
\end{quotation}
\vspace{6pt} The next result contains contains several easy-to-check
conditions that imply hypothesis \eqref{bite}.
\begin{proposition}[Sufficient conditions] \label{bitten}
Assume $\psi \in L^p, p \in (0,1)$, and that one of the following
conditions holds:
\begin{itemize}
\item[(a)] $\psi \in L^1$ with $\int_\Rd \psi \, dx \neq 0$;

\item[(b)] $\psi \in L^1$ with $\int_\Rd \psi \, dx = 0$, and $0
\not \equiv P\psi$ is real-valued and bounded either above or
below;

\item[(c)] $P\psi \in L^2_{loc}$ with $p \int_\cube |P\psi(b \, \cdot)|^2 \,
dx < (2-p) \left| \int_\cube (P\psi(b \, \cdot))^2 \, dx \right|$.
\end{itemize}

Then for some $\lambda \in \C$ (in fact with $|\lambda| < 1$),
\[
\| \lambda P\psi(b \, \cdot) - 1 \|_{L^p(\cube)} < 1 .
\]
This conclusion holds also if $p=1$ and $\psi$ satisfies condition
(a).
\end{proposition}
Section~\ref{bitten_proof} has the proof. Part (a) of the
proposition is essentially due to Filippov and Oswald
\cite[Lemma~1]{FO95}, and so is part (c) when $\psi$ is real-valued.
Part (b) is new.
\begin{quotation}
\noindent {\small \textsc{Examples.} Condition (c) in the
proposition holds if $\psi \in L^2$ is real valued with compact
support (or with rapid enough decay to ensure $P\psi \in L^2_{loc}$)
and with $P\psi \not \equiv 0$, because $p<2-p$. In particular, in
one dimension with $a_j=2^j$ and $b=1$, condition (c) (and also
condition (b)) covers the Haar wavelet
$\psi=\charfn_{[0,1/2)}-\charfn_{[1/2,1)}$. Thus while the dyadic
Haar system needs all scales $j \in \Z$ to span $L^2(\R)$, the small
scales $j>0$ suffice to span $L^p(\R)$ for $p \in (0,1)$, by
Theorem~\ref{surjectivity}. }
\end{quotation}

This Haar example reminds us the integral of $\psi$ can equal $0$,
in parts (b) and (c). When $p=1$, on the other hand, it is necessary
for Proposition~\ref{bitten} that $\psi$ have nonzero integral,
because otherwise
\[
\| \lambda P\psi(b \, \cdot) - 1 \|_{L^1(\cube)} \geq \left|
\int_\cube \left( \lambda P\psi(b \, \cdot) - 1 \right) dx \right| =
1 .
\]
Theorem~\ref{surjectivity} also requires $\psi$ to have nonzero
integral when $p=1$, because otherwise $S$ can synthesize only the
$L^1$ functions that have integral zero.

The next example shows it can be wise to ignore some of the
translations.
\begin{quotation}
\noindent {\small \textsc{Under-synthesizing example.} Hypothesis
\eqref{bite} definitely fails if the periodization of $\psi$
vanishes identically. But if the periodization of $\psi$ with
respect to some integer \emph{multiple} of $b$ is nontrivial, then
surjectivity of synthesis can still hold, as we explain. For
simplicity, work in one dimension with $b=1$ and suppose
$\psi=\eta-\eta(\cdot-1)$ for some $\eta \in L^p[0,1]$, so that
$\psi$ has the form of a unit step difference. (An illustrative
example is the Haar-type function
$\psi=\charfn_{[0,1)}-\charfn_{[1,2)}$.) The periodization of $\psi$
is $\eta-\eta \equiv 0$, and so Proposition~\ref{bitten} does not
apply. But if we consider the same $\psi$ with $b=2$, then the
$2$-periodization of $\psi$ equals $2\psi$ on the interval $[0,2)$.
Hence conditions (b) or (c) in Proposition~\ref{bitten} apply for $p
\in (0,1)$, assuming $\psi \not \equiv 0$ is either bounded or is
real valued and square integrable, respectively.
Theorem~\ref{surjectivity} with $b=2$ then tells us that every $f
\in L^p(\R)$ can be written as
\[
f = \sum_{j>0} \sum_{k \in \Z} \tilde{c}_{j,k} |a_j|^{1/p} \psi(a_j
x-2k)
\]
for some $\tilde{c} \in \ell^p(\Z_+ \times \Z)$. That is, $f=Sc$
where $c_{j,k}=\tilde{c}_{j,k/2}$ if $k$ is even, and $c_{j,k}=0$ if
$k$ is odd, and where $S$ denotes synthesis with $b=1$. In other
words, we have shown $S$ is surjective by under-synthesizing by a
factor of $2$, using only the even translates.}
\end{quotation}

Under-synthesis leads also to the following variant of
Theorem~\ref{surjectivity} and Proposition~\ref{bitten}, proved in
Section~\ref{undersynth_proof}.
\begin{theorem}[More synthesis onto $L^p$] \label{undersynth}
Assume $\psi \in L^p \cap L^1 \setminus \{ 0 \}$ for some $p \in
(0,1)$, and that $\psi$ is real-valued and has negative part
$\psi_-$ that is bounded with compact support.

Then $S : \ell^p(\Z_+ \times \Zd) \to L^p$ is open, and surjective.
Indeed, a constant $C=C(\psi,b,p)$ exists such that if $f \in L^p$
then there is a sequence $c \in \ell^p(\Z_+ \times \Zd)$ with $Sc=f$
and $\| c \|_{\ell^p(\Z_+\times \Zd)} \leq C \| f \|_p$.
\end{theorem}
The conclusion holds also if the positive part of $\psi$ is bounded
and has compact support, instead of the negative part.
\begin{quotation}
\noindent {\small \textsc{Example.} If $\eta \in L^\infty[0,1]$ is
real-valued and $\eta \not \equiv 0$, then the second difference
function $\psi(x)=\eta(x+1)-2\eta(x)+\eta(x-1)$ is covered by
Theorem~\ref{undersynth}, but $\psi$ is not covered by
Theorem~\ref{surjectivity} when $b=1$ because $P\psi \equiv 0$.}
\end{quotation}

\vspace{6pt} Equivalence of the $L^p$ and $\ell^p$ metrics follows
immediately from Proposition~\ref{continuitysynth} and
Theorem~\ref{surjectivity}:
\begin{corollary}[Affine atomic decomposition of $L^p$] \label{atomic}
Assume $\psi \in L^p, p \in (0,1]$, and that $\psi$ satisfies
hypothesis \eqref{bite} in Theorem~\ref{surjectivity}. Then for all
$f \in L^p$,
\[
\| f \|_p \approx \inf \left\{ \| c \|_{\ell^p(\Z_+ \times \Zd)} :
f = Sc \right\} .
\]
\end{corollary}
For $p=1$, the corollary was proved by Bruna \cite[Theorem~4]{B06}
for $\psi \in L^1$ with $\int_\Rd \psi \, dx \neq 0$. His duality
methods apply without our assumption that the translations lie in a
lattice.

Next we localize the atomic decomposition to an open set $\Omega
\subset \Rd$.

\subsubsection*{Definition.} Say that a sequence $c =
\{ c_{j,k} \}_{j>0, k\in \Zd}$ is \emph{adapted to $\Omega$ and
$\psi$} if $\spt(\psi_{j,k}) \subset \Omega$ whenever $c_{j,k} \neq
0$, or in other words if $c_{j,k}=0$ whenever $\spt(\psi_{j,k}) \cap
\Omega^c \neq \emptyset$.

\vspace{6pt} The purpose of the definition is to ensure $Sc=0$ on
the complement of $\Omega$.
\begin{corollary}[Affine atomic decomposition of $L^p(\Omega)$] \label{atomicdomain}
Assume $\Omega \subset \Rd$ is open and nonempty, take $p \in
(0,1]$, and suppose $\psi \in L^p$ is compactly supported and
satisfies hypothesis \eqref{bite} in Theorem~\ref{surjectivity}.
Then for all $f \in L^p(\Omega)$,
\[
\| f \|_{L^p(\Omega)} \approx \inf \left\{ \| c \|_{\ell^p(\Z_+
\times \Zd)} : \text{$f = Sc$ and $c$ is adapted to $\Omega$ and
$\psi$} \right\} .
\]
\end{corollary}
The constants in this metric equivalence are the same as in
Corollary~\ref{atomic}; thus they depend on $\psi, b$ and $p$ but
are independent of $\Omega$. The corollary is proved in
Section~\ref{atomicdomain_proof}.

The \emph{full analysis operator} puts the analysis sequences from
all the different scales into a combined doubly-indexed sequence, by
\[
Tf = \{ |\det b| \Theta^{-1} \langle \Theta f , \phi(a_j \, \cdot -
bk) \rangle \}_{j>0, k \in \Zd} .
\]
It too yields a coefficient metric, as we prove in
Section~\ref{coeffanal_proof}:
\begin{corollary}[Analysis metric for $L^p$] \label{coeffanal}
Assume $p \in (0,1]$, and take an analyzer $\phi \in L^\infty$ with
$P|\phi| \in L^\infty$ and $\int_\Rd \phi \, dx = 1$.

Then for all $f \in L^p$,
\[
\|f\|_p \approx \| Tf \|_{\ell^\infty(\ell^p)} = \sup_{j>0} \, [
\sum_{k \in \Zd} | \langle \Theta f , \phi(a_j \cdot - bk) \rangle|
\, ]^{1/p} \, |\det b|.
\]
\end{corollary}

\vspace{6pt}
\subsubsection*{Synthesis at a single scale: conditions for a $p$-Riesz basis.}
This paper concentrates on synthesizing $L^p$ by using all scales
$j>0$. But we divert briefly now from this mission to give a fuller
account of the synthesis operator at a single scale $j$. We
establish conditions for $S_j$ to be injective and have closed
range, which is equivalent to ``stability'' or a $p$-Riesz basis
condition, or an atomic decomposition at scale $j$.

Write $e_\xi(x)=e^{2\pi i \xi x}$, where $\xi \in \Rd$ is a row and
$x \in \Rd$ is a column vector.
\begin{theorem}[$p$-Riesz basis at scale $j$] \label{inject}
Assume $\psi \in L^p, p \in (0,1]$, and that $P(e_\xi \psi) \not
\equiv 0$ for each $\xi \in \Rd$. Let $j>0$.

Then $S_j : \ell^p(\Zd) \to L^p$ is injective.

If in addition $\psi$ has compact support, then $\| S_j s \|_p
\approx \| s \|_{\ell^p(\Zd)}$, or
\begin{equation} \label{inject1}
\| \sum_{k \in \Zd} s_k \psi_{j,k} \|_p \approx ( \sum_{k \in \Zd}
|s_k|^p )^{1/p} , \qquad s \in \ell^p(\Zd) .
\end{equation}
That is, the collection of integer shifts $\{ \psi_{j,k} : k \in
\Zd \}$ forms a $p$-Riesz basis for its span. Hence the range
$S_j(\ell^p(\Zd))$ is closed in $L^p$.
\end{theorem}
The constants in the norm equivalence \eqref{inject1} depend on
$\psi, b$ and $p$, but are independent of $j$.

Theorem~\ref{inject} is proved in Section~\ref{inject_proof}. The
periodization hypothesis $P(e_\xi \psi) \not \equiv 0$ says in the
Fourier domain (when $\psi$ has a Fourier transform) that the
sequence $\{ \widehat{\psi}(\ell b^{-1}-\xi) \}_{\ell \in \Zd}$ is
nontrivial, because this sequence gives the Fourier coefficients of
the $\Zd$-periodic function $P(e_\xi \psi)(bx)$.

Ron \cite{R89} proved the injectivity conclusion in
Theorem~\ref{inject}, and Jia \cite{J98} proved the Riesz basis
conclusion \eqref{inject1}, except that here we work directly with
the periodization hypothesis $P(e_\xi \psi) \not \equiv 0$ and thus
need not assume like Ron and Jia that $\psi$ can be Fourier
transformed. (Note our $\psi \in L^p$ need not be locally
integrable, or even a distribution.) Otherwise we follow Ron and
Jia's method exactly. Incidentally, Ron obtained a converse when
$\psi$ has compact support and the restriction that $s \in
\ell^p(\Zd)$ is dropped, saying that injectivity implies nontrivial
periodizations, and Jia proved a converse saying that if
\eqref{inject1} holds then $P(e_\xi \psi) \not \equiv 0$ for all
$\xi \in \Rd$.

\subsubsection*{Properties of the class of
synthesizers.} So far we have concentrated on individual
synthesizers and analyzers. Now we broaden the view and consider the
whole class
\[
{\mathcal S}^p = \{ \psi \in L^p : \text{$S_\psi$ maps $\ell^p(\Z_+
\times \Zd)$ onto $L^p$} \}
\]
of surjective synthesizers, where the $\psi$-dependence of the
synthesis operator is emphasized by writing $S=S_\psi$.
\begin{theorem}[Most $L^p$ functions are surjective
affine synthesizers] \label{properties} Let $p \in (0,1]$. Then
${\mathcal S}^p$ is dense, open, and path connected in $L^p$.
\end{theorem}
The proof is in Section~\ref{properties_proof}.

\section{\bf Remarks on $L^p$ for $p \geq 1$, and on Hardy and Sobolev
spaces.} \label{remarksLHS}

This paper treats affine synthesis in $L^p$ for $p \in (0,1]$. It is
helpful to contrast these results and methods with the corresponding
work for $p \geq 1$.

Affine synthesis when $p \geq 1$ was treated in my paper \cite{bl5}
with H.-Q. Bui. The domain of synthesis there is a mixed-norm
sequence space $\ell^1(\ell^p)$ (meaning $\ell^p$ with respect to
translations $k$ and then $\ell^1$ with respect to dilations $j$).
The only situations I know where $p
> 1$ and synthesis is bounded on the domain $\ell^p(\Z_+ \times
\Zd)$ like in this paper (meaning $\ell^p$ with respect to both
translation and dilation) are when the $\psi_{j,k}$ possess some
``cancelation between scales'', such as for wavelets in $L^2$ where
the $\psi_{j,k}$ form an orthonormal basis. This is discussed in
\cite[\S8.1]{bl5}. Of course, the two domain spaces coincide when
$p=1$.

Surjectivity of synthesis was established in \cite{bl5} via
$L^p$-controlled approximation, using linear analysis and linear
synthesis. (A linear analysis operator could be used because $L^p$
functions are locally integrable when $p \geq 1$.) But rather than
obtaining surjectivity of synthesis from Theorem~\ref{sample} and
Proposition~\ref{bitten}(a) (which both still hold for $p \geq 1$
assuming $P|\psi| \in L^p_{loc}$), we used in \cite{bl5} our
scale-averaging technique from \cite{bl2} to essentially take
$\sigma=0$, thereby improving the estimate on the norm of $c$ in the
analogue of Theorem~\ref{surjectivity}. Scale-averaging of this kind
completely fails when $p \in (0,1)$, because the unit ball of $L^p$
is non-convex. Hence we must make do in this paper with the somewhat
cruder approximation provided by Theorem~\ref{sample}.

Surjectivity of synthesis when $\psi$ has nonzero integral was
proved earlier in \cite[Theorem~4]{B06} for $p=1$, and even earlier
in \cite[Theorem~2]{T99},\cite{T99b} for $p \geq 1$. These authors
both proceed by studying the analysis operator (proving $\| Tf
\|_{\ell^\infty(\ell^p)} \approx \| f \|_p$) and then invoking
duality; thus they provide no constructive method of synthesis like
we provide in \cite{bl5} and in this paper.

The theory of $p$-Riesz bases for $1<p<\infty$ is described in
Christensen's book \cite[\S7.2,17.4]{Ch03}, and Theorem~\ref{inject}
should be read in that context.

A density result for synthesizers is known when $p
> 1$, like in Theorem~\ref{properties}. But openness and path connectedness are
not known: see \cite[\S4]{bl2} for relevant remarks. Incidentally,
some path connectedness results in the difficult wavelet and wavelet
frame case, in $L^2$, can be found in \cite{B05,W98} and the
references therein.

Scale-averaging and surjectivity results for the Hardy space $H^1$
and for Sobolev spaces are contained in \cite{bl3,bl4,bl5}. Many
smooth wavelet systems are known to span the Hardy space too
\cite[Theorem~5.6.19]{HW96}. But note the Haar system does not span
Hardy space: the closed span of the Haar system in $H^1(\R)$ is the
proper subspace $\{ f \in H^1(\R) : \int_0^\infty f \, dx = 0 \}$,
by \cite[Theorem~2.1]{AST06}.

\section{\bf Open problems}
\label{openproblems}

\subsection*{Obtain a larger class of synthesizers?}
The synthesizer $\psi \in L^p$ in Theorem~\ref{surjectivity} is
assumed to satisfy
\begin{equation} \label{biteagain}
\| \lambda P\psi(b \, \cdot) - 1 \|_{L^p(\cube)} < 1 \qquad
\text{for some $\lambda \in \C$.}
\end{equation}
Can this condition be weakened (or even eliminated), to obtain a
larger class of synthesizers?

Two facts advocate for caution here. First, our sufficient condition
\eqref{biteagain} fails for some $\psi \in L^p$, as Tachev
\cite{T95} observed in one dimension (see Section~\ref{tachev_proof}
below). Second, if $\int_\Rd \psi \, dx = 0$ and $P\psi \in
L^\infty$ satisfies the reverse of the inequality in
Proposition~\ref{bitten}(c), then the proof of that part of the
proposition adapts to give the reverse of \eqref{biteagain} for all
small $\lambda \in \C, \lambda \neq 0$.

\subsection*{Employ large scales?}
We have only used small scales $j > 0$ in this paper. The large
scales $j \leq 0$ are employed (for example) by wavelet systems in
$L^p, p>1$, and so it is natural to ask how our synthesis results
for $p \in (0,1)$  might be affected by the introduction of large
scales.

For the sake of concreteness, consider dyadic dilations $a_j = 2^j
I$ for $j \in \Z$. Then the synthesis operator $S : \ell^p(\Z \times
\Zd) \to L^p$ is bounded, since it makes no difference in the proof
of Proposition~\ref{continuitysynth} that the dilation parameter
runs over all $j \in \Z$. The hope is to find a new approach to
proving surjectivity of $S$ that somehow incorporates the large
scales, and thereby yields synthesizers not already covered by
Proposition~\ref{bitten} and Theorem~\ref{undersynth}.

In addition to introducing large scales, one might allow modulations
(wavepacket theory, like in \cite{CKS06,FF06,LWW04}) or
multiplication by polynomials (like in Gausslet theory \cite{T02}),
and again ask for synthesizers that are surjective onto $L^p$ for $p
\in (0,1)$.

\subsection*{Establish a rate of convergence?}
Theorem~\ref{sample} says that if $\psi$ has constant periodization
$P \psi \equiv 1$, then nonlinear quasi-interpolation at small
scales will converge to the sampled function, meaning $S_j T_j f \to
f$ in $L^p$.

How fast is this convergence to $f$? When $p \geq 1$, Strang--Fix
theory \cite{BJ85,bl4,HR05,SF73} provides very precise convergence
rates, for synthesizers $\psi$ satisfying Strang--Fix conditions and
signals $f$ belonging to a Sobolev space. The challenge is to
develop similar approximation rate results when $p \in (0,1)$,
bearing in mind the nonlinearity of the analysis operator.

\subsection*{Mexican hat spanning problem}
It is an open problem of Y. Meyer \cite[p.~137]{M92} to determine
whether the affine system $\{ \psi(2^j x-k) : j, k \in \Z \}$ spans
$L^p(\R)$ for each $p \in (1,\infty)$, when
$\psi(x)=(1-x^2)e^{-x^2/2}$ is the Mexican hat function (the second
derivative of the Gaussian). This is known to be true in $L^2(\R)$,
where the system forms a frame, but remains open for all other
$p$-values between $1$ and $\infty$.

Theorem~\ref{surjectivity} completely resolves the Mexican hat
problem for $0<p<1$, because the Mexican hat $\psi$ belongs to $L^p$
and has periodization $P\psi \in L^\infty \setminus \{ 0 \}$, so
that it satisfies the hypotheses of Proposition~\ref{bitten}(b) and
hence of Theorem~\ref{surjectivity}. (To see $P\psi \not \equiv 0$,
note the Fourier coefficients of $P\psi$ are given by the values of
$\widehat{\psi}$ at the integers, which are nonzero because
$\widehat{\psi}(\xi)$ equals a constant times $\xi^2$ times a
Gaussian.) In fact, Theorem~\ref{surjectivity} shows the Mexican hat
system spans $L^p$ using only the small scales $j>0$, rather than
all $j \in \Z$ as allowed by the problem.

More discussion of the Mexican hat problem can be found in
\cite[\S7]{bl5}.

\section{\bf Radial stretching}
\label{elementary}

We will need some elementary properties of the radial stretch
function $\Theta(z) = |z|^{p-1} z$. Remember $|\Theta (z)|=|z|^p$.
\begin{lemma}[H\"{o}lder continuity] \label{holder}
If $p \in (0,1]$ then
\[
|\Theta(w) - \Theta(z)| \leq 2^{1-p} |w-z|^p , \qquad w,z \in \C.
\]
\end{lemma}
\begin{proof}[Proof of Lemma~\ref{holder}]
First reduce to the case $w=1$ and $z=re^{i\phi}$ with $0 < r \leq
1$ and $0 \leq \phi \leq \pi$, by a rotation, dilation and
reflection if necessary. Thus we want
\begin{equation} \label{holder1}
\left| \frac{1-r^p e^{i\phi}}{2} \right| \leq \left| \frac{1-r
e^{i\phi}}{2} \right|^p .
\end{equation}
This holds for $\phi=\pi$ just by concavity of $r \mapsto r^p$. To
prove \eqref{holder1} for $0 \leq \phi < \pi$, we write $R(\phi)$
for the ratio of the left side of \eqref{holder1} over the right
side, and observe that
\begin{align*}
\frac{\partial \ }{\partial \phi} \log R(\phi) & = \sin \phi \left\{
r^p |1-r^p e^{i\phi}|^{-2} - pr|1-r e^{i\phi}|^{-2} \right\} \\
& \geq \sin \phi \left\{
r^p |1-r^p e^{i\phi}|^{-2} - r|1-r e^{i\phi}|^{-2} \right\} \qquad \text{since $p \leq 1$} \\
& = \sin \phi \left\{ (r^{-1}+r) - (r^{-p}+r^p) \right\} r^p |1-r^p
e^{i\phi}|^{-2} r|1-r e^{i\phi}|^{-2} \\
& \geq 0
\end{align*}
since $0 < r \leq r^p \leq 1$.
\end{proof}
\begin{lemma}[Local Lipschitz continuity of the inverse] \label{lipschitz}
If $p \in (0,1]$ then
\[
|\Theta^{-1}(w) - \Theta^{-1}(z)| \leq \frac{1}{p} |w - z| \,
\max(|w|,|z|)^{(1-p)/p} , \qquad w,z \in \C.
\]
\end{lemma}
\begin{proof}[Proof of Lemma~\ref{lipschitz}]
Again we can reduce to the case $w=1$ and $z=re^{i\phi}$ with $0
\leq r \leq 1$ and $0 \leq \phi \leq \pi$. Thus we want
\[
|1-r^{1/p} e^{i\phi}| \leq \frac{1}{p} |1-r e^{i\phi}| .
\]
This holds at $\phi=0$ by the mean value theorem, since $(1/p) \geq
1$. To prove it for $0 < \phi \leq \pi$, we compute that
\[
\frac{\partial \ }{\partial \phi} |1-r^{1/p} e^{i\phi}|^2 = 2r^{1/p}
\sin \phi \leq 2p^{-2} r \sin \phi = \frac{\partial \ }{\partial
\phi} \frac{1}{p^2} |1-r e^{i\phi}|^2,
\]
since $p \leq 1$ and $0 \leq r \leq 1$.
\end{proof}

\section{\bf $L^p$ proofs}
\label{lebesgue_proofs}

\subsection{\bf Proof of Proposition~\ref{continuitysynth} ---
$\text{synthesis\ } \ell^p \to L^p$} \
\label{continuitysynth_proof}

We have
\begin{align*}
\| Sc \|_p^p & \leq \int_\Rd ( \sum_{j>0} \sum_{k \in \Zd}
|c_{j,k} \psi_{j,k}(x)| )^p \, dx \\
& \leq \int_\Rd \sum_{j>0} \sum_{k \in \Zd} |c_{j,k}
\psi_{j,k}(x)|^p \, dx \qquad \text{since $p \in (0,1]$} \\
& = \sum_{j>0} \sum_{k \in \Zd} |c_{j,k}|^p \| \psi_{j,k} \|_p^p =
\| c \|_{\ell^p(\Z_+ \times \Zd)}^p \| \psi \|_p^p
\end{align*}
since $\| \psi_{j,k} \|_p = \| \psi \|_p$ for all $j,k$ by our
normalization of $\psi_{j,k}$. Thus the series for $Sc$ converges
a.e.\ to an $L^p$-function, which implies unconditional convergence
of the series in $L^p$ (with the help of the dominated convergence
theorem).

\subsubsection*{Aside.} Obviously all one really needs
here, in order for the synthesis operator to be continuous, is that
the synthesizing collection $\{ \psi_{j,k} \}$ be bounded in $L^p$.

\subsection{\bf Proof of Theorem~\ref{continuityanal} ---
$\text{analysis\ } L^p \to \ell^p$} \ \label{continuityanal_proof}

Take $f \in L^p$ and define a sequence
\[
s_k = \langle \Theta f , \phi(a_j \, \cdot - bk) \rangle , \qquad k
\in \Zd .
\]
Then
\begin{align*}
|\det b| \sum_{k \in \Zd} |s_k| & \leq \int_\Rd |\Theta f(y)| |\det b| \sum_{k \in \Zd} |\phi(a_j y - bk)| \, dy \\
& \leq \| f \|_p^p \| P|\phi| \|_\infty .
\end{align*}
Since $(T_j f)_k = |\det b| \Theta^{-1}(s_k)$ by definition, we
deduce $T_j$ maps $L^p$ into $\ell^p(\Zd)$ with
\begin{align}
\| T_j f \|_{\ell^p(\Zd)} & = |\det b| ( \sum_{k \in \Zd}
|s_k| )^{1/p} \notag \\
& \leq |\det b|^{1-1/p} \| P|\phi| \|_\infty^{1/p} \| f \|_p .
\label{coeffexplicit}
\end{align}

To obtain H\"{o}lder continuity of the analysis operator, we
consider $g \in L^p$ and let $t_k = \langle \Theta g , \phi(a_j \,
\cdot - bk) \rangle$. Then
\begin{align}
d_{\ell^p}(T_j f , T_j g)
& = \| T_j f - T_j g \|_{\ell^p(\Zd)}^p \notag \\
& = \sum_{k \in \Zd} |\det b|^p |\Theta^{-1}(s_k) - \Theta^{-1}(t_k)|^p \notag \\
& \leq p^{-p} \sum_{k \in \Zd} |\det b|^p |s_k - t_k|^p
\max(|s_k|,|t_k|)^{1-p} \qquad \text{by Lemma~\ref{lipschitz}} \notag \\
& \leq p^{-p} ( \sum_{k \in \Zd} |\det b| |s_k - t_k| )^p ( \sum_{k
\in \Zd} \max(|s_k|,|t_k|) )^{1-p} \label{holdereq1}
\end{align}
by H\"{o}lder's inequality on the sum. The first sum is bounded by
\begin{align*}
|\det b| \sum_{k \in \Zd} |s_k - t_k| & \leq |\det b| \sum_{k \in
\Zd} \int_\Rd |\Theta f(x) - \Theta g(x)| |\phi(a_j x - bk)| \, dx \\
& \leq 2^{1-p} \int_\Rd |f(x)-g(x)|^p \, dx \, \| P|\phi| \|_\infty
\qquad \text{by Lemma~\ref{holder}} \\
& = 2^{1-p} \| P|\phi| \|_\infty d_p(f,g) .
\end{align*}
And the second sum is bounded by
\[
\sum_{k \in \Zd} (|s_k|+|t_k|) \leq |\det b|^{-1} \| P|\phi|
\|_\infty [d_p(f,0) + d_p(0,g)] .
\]
Combining these estimates proves the desired H\"{o}lder continuity
in the theorem, with $C=2^{p(1-p)} p^{-p} |\det b|^{p-1} \| P|\phi|
\|_\infty$.

\subsection{\bf Continuity of the analysis operator, with respect to the analyzer} \ \label{contwrtanal_sec}

The preceding section proved continuity of the analysis operator
with respect to the signal $f$. Now we show it is continuous with
respect to the analyzer $\phi$. Both results will be used in the
proof of Theorem~\ref{sample}.

To emphasize that the analysis operator depends on $\phi$, we write
$T_j = T_{j,\phi}$ in this section.
\begin{lemma} \label{contwrtanal}
Assume $p \in (0,1]$, and take $f \in L^\infty$ with compact support
in $\Rd$.

Then for each $j$, the map $\phi \mapsto T_{j,\phi} f$ is locally
H\"{o}lder continuous from $L^1$ to $\ell^p(\Zd)$, with
\[
d_{\ell^p}(T_{j,\phi} f , T_{j,\varphi} f) \leq C |\det b|^{p-1} \|
f \|_\infty^p (\| \phi \|_1 + \| \varphi \|_1 )^{1-p} \| \phi -
\varphi \|_1^p , \qquad \phi, \varphi \in L^1 .
\]
Here $C$ depends on the exponent $p$, the support of $f$ and on
$\max_{J>0} \| a_J^{-1} \|$, but not on the dilation scale $j$.
\end{lemma}
\begin{proof}[Proof of Lemma~\ref{contwrtanal}]
First we show the analysis operator is well defined under the
hypotheses of this lemma. Again write $s_k = \langle \Theta f ,
\phi(a_j \, \cdot - bk) \rangle$, so that $(T_{j,\phi} f)_k = |\det
b| \Theta^{-1}(s_k)$. Then
\begin{align}
|\det b| \sum_{k \in \Zd} |s_k|
& = \int_\Rd \left[ |\det a_j^{-1}b| \sum_{k \in \Zd} |f(a_j^{-1}y+a_j^{-1}bk)|^p \right] |\phi(y)| \, dy \notag \\
& \qquad \qquad \qquad \qquad \text{by $y \mapsto a_j^{-1}(y+bk)$ in the integral for $s_k$} \notag \\
& \leq C \| f \|_\infty^p \| \phi \|_1 , \label{wrt1}
\end{align}
where the constant $C$ comes from estimating the Riemann sum and
thus depends only on the diameter of the support of $f$ and on the
``step size'' $\| a_j^{-1} \|$. Hence $T_{j,\phi} f$ is well defined
and belongs to $\ell^p(\Zd)$, with
\begin{align*}
\| T_{j,\phi} f \|_{\ell^p(\Zd)} & = |\det b| ( \sum_{k \in \Zd} |
s_k | )^{1/p} \\
& \leq C |\det b|^{1-1/p} \| f \|_\infty \| \phi \|_1^{1/p}
\end{align*}
holding whenever $\phi \in L^1$ and $f \in L^\infty$ has compact
support.

To obtain H\"{o}lder continuity we now consider another analyzer
$\varphi \in L^1$ and define the corresponding sequence $t_k =
\langle \Theta f , \varphi(a_j \, \cdot - bk) \rangle$, so that
$(T_{j,\varphi} f)_k = |\det b| \Theta^{-1}(t_k)$. Then
\[
d_{\ell^p}(T_{j,\phi} f , T_{j,\varphi} f) \leq p^{-p} ( \sum_{k \in
\Zd} |\det b| |s_k - t_k| )^p ( \sum_{k \in \Zd} \max(|s_k|,|t_k|)
)^{1-p}
\]
by arguing like for \eqref{holdereq1} in the proof of
Theorem~\ref{continuityanal}. The first sum is bounded by
\[
|\det b| \sum_{k \in \Zd} |s_k - t_k| \leq C \| f \|_\infty^p \|
\phi - \varphi \|_1 ,
\]
as one sees by applying estimate \eqref{wrt1} to $\phi-\varphi$
instead of to $\phi$. The second sum is bounded by
\[
\sum_{k \in \Zd} (|s_k|+|t_k|) \leq C |\det b|^{-1} \| f \|_\infty^p
( \| \phi \|_1 + \| \varphi \|_1 ) ,
\]
by \eqref{wrt1}. Combining these estimates proves the desired
H\"{o}lder continuity in the lemma.
\end{proof}

\subsection{\bf Proof of Theorem~\ref{sample} --- affine quasi-interpolation}
\ \label{sample_proof}

[Step 1.] Let $f \in L^p$. First we reduce to $f$ being continuous
with compact support. Indeed, given $\e>0$ we can choose a
continuous function $g$ with compact support and $d_p(f,g) < \e$.
Then
\begin{align*}
d_p(S_j T_j f, S_j T_j g) & \leq \| \psi \|_p^p
d_{\ell^p}(T_j f, T_j g) && \text{by Proposition~\ref{continuitysynth}} \\
& \leq C(b,p,\phi,\psi) [d_p(f,0) + d_p(0,g)]^{1-p} d_p(f,g)^p
&& \text{by Theorem~\ref{continuityanal}} \\
& \leq C(b,p,\phi,\psi) [2d_p(f,0) + \e]^{1-p} \e^p .
\end{align*}
So if we prove $\lim_{j \to \infty} \| S_j T_j g - g \|_p^p = \sigma
\| g \|_p^p$ where
\[
\sigma = \| P\psi(b \, \cdot) - 1 \|_{L^p(\cube)}^p,
\]
then it follows that $\lim_{j \to \infty} \| S_j T_j f - f \|_p^p =
\sigma \| f \|_p^p$ as desired, by taking $\e$ arbitrarily small.
Thus we may assume from now on that $f$ is  continuous with compact
support.

\vspace{6pt} [Step 2.] Now we reduce to $\phi \in L^\infty$ having
compact support. The analyzer $\phi$ is certainly integrable, since
the periodization $P|\phi|$ is assumed to be bounded and hence is
locally integrable. Therefore given $\e>0$ we can choose $\varphi
\in L^\infty$ with compact support and $\| \phi - \varphi \|_1 < \e$
and $\int_\Rd \varphi \, dx = 1$. Proceeding analogously to the
reduction in Step~1, we observe
\begin{align*}
d_p(S_j T_{j,\phi} f, S_j T_{j,\varphi} f) & \leq \| \psi \|_p^p
d_{\ell^p}(T_{j,\phi} f, T_{j,\varphi} f) && \text{by Proposition~\ref{continuitysynth}} \\
& \leq C(b,p,\psi,f) ( \| \phi \|_1 + \| \varphi \|_1 )^{1-p} \|
\phi - \varphi \|_1^p
&& \text{by Lemma~\ref{contwrtanal}} \\
& \leq C(b,p,\psi,f) (2\| \phi \|_1 + \e)^{1-p} \e^p .
\end{align*}
Thus we need only prove $\lim_{j \to \infty} \| S_j T_{j,\varphi} f
- f \|_p^p = \sigma \| f \|_p^p$, because then taking $\e$
arbitrarily small implies the corresponding limit with $\phi$
instead of $\varphi$. Thus we may assume from now on that $\phi \in
L^\infty$ has compact support.

\vspace{6pt} [Step 3.] Next we reduce to analyzing $f$ with
pointwise sampling. Begin by uniformly sampling the continuous
function $f$ at scale $j$, and recording the results in the sequence
\[
U_j f = \{ |\det a_j|^{-1/p} |\det b| f(a_j^{-1}bk) \}_{k \in \Zd} .
\]
That is, $U_j$ is a \emph{pointwise analysis operator at scale $j$}.
We aim to show average sampling and pointwise sampling are the same
in the limit, or more precisely that
\begin{equation} \label{TUclose}
\| T_j f - U_j f \|_{\ell^p(\Zd)} \to 0 \qquad \text{as $j \to
\infty$.}
\end{equation}

Take $j$ large enough that
\begin{equation} \label{supportest}
|a_j^{-1} y| < 1 \qquad \text{for all $y \in \spt(\phi)$,}
\end{equation}
using here that $\phi$ is compactly supported and $\| a_j^{-1} \|
\to 0$. Write $F_r$ for the set of points within distance $r>0$ of
the support of $f$, and let
\[
K(j) = \{ k \in \Zd : a_j^{-1}bk \in F_1 \} .
\]
Then
\begin{equation} \label{nearby1}
0 \neq (U_j f)_k \quad \Longrightarrow \quad k \in K(j) ,
\end{equation}
because if $f(a_j^{-1}bk) \neq 0$ then $a_j^{-1}bk \in \spt(f)
\subset F_1$. And
\begin{equation} \label{nearby2}
0 \neq (T_j f)_k = |\det b| \Theta^{-1} \langle \Theta f , \phi(a_j
\, \cdot - bk) \rangle \quad \Longrightarrow \quad k \in K(j) ,
\end{equation}
because if $(T_j f)_k \neq 0$ then there exists $x \in \spt(f)$ with
$a_j x - bk \in \spt(\phi)$, so that $|x-a_j^{-1}bk|<1$ by
\eqref{supportest}, which implies $a_j^{-1}bk \in F_1$ and hence $k
\in K(j)$.

In view of \eqref{nearby1} and \eqref{nearby2}, when proving
\eqref{TUclose} we need only sum over $K(j)$. Thus
\begin{align*}
& \| T_j f - U_j f \|_{\ell^p(\Zd)}^p \\
& = |\det b|^{p-1} |\det a_j^{-1} b| \sum_{k \in K(j)} |
\Theta^{-1} (|\det a_j| \langle
\Theta f , \phi(a_j \, \cdot - bk) \rangle) - f(a_j^{-1}bk)|^p \\
& \leq |\det b|^{p-1} |\det a_j^{-1} b| \# K(j) \cdot M(j)^p
\end{align*}
where
\[
M(j) = \sup_{k \in \Zd} |\Theta^{-1}(\int_\Rd \Theta
f(a_j^{-1}(y+bk)) \overline{\phi(y)} \, dy) - \Theta^{-1}(\Theta
f(a_j^{-1}bk))| .
\]
Since $|\det a_j^{-1}b| \# K(j)$ is bounded by the volume of
$F_2$, for all large $j$, we can see that in order to prove
\eqref{TUclose} it suffices to show $M(j) \to 0$.

Notice the arguments of $\Theta^{-1}(\ldots)$ in the definition of
$M(j)$ are bounded independently of $j$ and $k$, since $f$ is
bounded and $\phi$ is integrable. Hence the convergence of $M(j)$ to
$0$ follows from local uniform continuity of $\Theta^{-1}$, since
the distance between the arguments converges to $0$ uniformly with
respect to $k$, as follows:
\begin{align*}
& \sup_{k \in \Zd} | \int_\Rd \Theta f(a_j^{-1}(y+bk))
\overline{\phi(y)} \, dy - \Theta f(a_j^{-1}bk)| \\
& \leq \int_\Rd \sup_{k \in \Zd} |\Theta f(a_j^{-1}y+a_j^{-1}bk) -
\Theta f(a_j^{-1}bk)| |\phi(y)| \, dy
\qquad \text{since $\int_\Rd \phi \, dx = 1$} \\
& \to 0
\end{align*}
as $j \to \infty$, with dominated convergence justified by uniform
continuity of the compactly supported function $\Theta f$ and
integrability of $\phi$. This proves $M(j) \to 0$, and hence proves
\eqref{TUclose}.

\vspace{6pt} [Step 4.] We next derive the theorem with $U_j$ in
place of $T_j$, in other words, we show
\begin{equation} \label{pointversion}
\lim_{j \to \infty} \| S_j U_j f - f \|_p^p = \sigma \| f \|_p^p .
\end{equation}
This implies the theorem because
\begin{align*}
\| S_j T_j f - S_j U_j f \|_p & \leq \| \psi \|_p \| T_j f - U_j f
\|_{\ell^p(\Zd)}
&& \text{by Proposition~\ref{continuitysynth}} \\
& \to 0 && \text{by \eqref{TUclose}.}
\end{align*}

To prove \eqref{pointversion}, we start by decomposing
\begin{equation} \label{decomp}
(S_j U_j f)(x) - f(x) = [P\psi(a_j x) - 1] f(x) + \text{Rem}_j(x) ,
\end{equation}
where the remainder is
\[
\text{Rem}_j(x) = |\det b| \sum_{k \in \Zd} ( f(a_j^{-1}bk) - f(x)
) \psi(a_j x - bk) .
\]
The first term on the right of \eqref{decomp} has limit
\[
\lim_{j \to \infty} \| (P\psi(a_j \, \cdot) - 1) f \|_p^p =
(\text{mean value of $|P\psi - 1|^p$}) \cdot \| f \|_p^p = \sigma \|
f \|_p^p
\]
by a Riemann--Lebesgue argument, since $|P\psi(a_j \cdot) - 1|^p$
oscillates rapidly around its mean value when $j$ is large. For
details see \cite[Lemma~26]{bl1}, for example, noting that $|P\psi -
1|^p$ is locally integrable and $|f|^p$ is bounded with compact
support.

\vspace{6pt} [Step 5.] To complete the proof of \eqref{pointversion}
we have only to show the remainder term $\text{Rem}_j$  vanishes in
the limit, in $L^p$. We have
\[
|\text{Rem}_j(x)|^p \leq |\det b|^p \sum_{k \in \Zd} |
f(a_j^{-1}bk) - f(x) |^p |\psi(a_j x - bk)|^p ,
\]
and so after integrating with respect to $x$ and making the change
of variable $x \mapsto a_j^{-1}(x+bk)$, we find
\[
\| \text{Rem}_j \|_p^p \leq |\det b|^{p-1} \int_\Rd R_j(x)
|\psi(x)|^p \, dx
\]
where
\[
R_j(x) = |\det a_j^{-1}b| \sum_{k \in \Zd} | f(a_j^{-1}bk) -
f(a_j^{-1}(x+bk)) |^p .
\]

Formally, $\Rem_j \to 0$ in $L^p$ because $R_j(x)$ is a Riemann sum
that passes in the limit to the integral
\[
\int_\Rd | f(z) - f(0+z) |^p dz = 0 .
\]
To prove $\Rem_j \to 0$ rigorously by dominated convergence, it is
enough to show $R_j(x) \to 0$ pointwise and that $R_j$ is bounded
independently of $x$ and $j$, for all large $j$. To get boundedness
of $R_j$, we estimate that
\begin{align*}
|R_j(x)| & \leq 2 \max_{z \in \Rd} |\det a_j^{-1}b| \sum_{k
\in \Zd} |f(z+a_j^{-1}bk)|^p \qquad \text{for all $x \in \Rd$} \\
& \to 2 \| f \|_p^p
\end{align*}
as $j \to \infty$, by a Riemann sum argument applied to the
continuous,  compactly supported function $|f|^p$. Thus $R_j$ is
bounded independently of $x$ and $j$, for all large $j$.

To get that $R_j(x) \to 0$ pointwise, we fix $x \in \Rd$ for the
rest of the proof, and take $j$ to be large enough that $|a_j^{-1}
x| < 1$. Then we need only sum over $K(j)$, when we evaluate
$R_j(x)$:
\begin{align*}
R_j(x) & = |\det a_j^{-1}b| \sum_{k \in K(j)} | f(a_j^{-1}bk) -
f(a_j^{-1}(x+bk)) |^p \\
& \leq |\det a_j^{-1} b| \# K(j) \cdot N(j)^p
\end{align*}
where
\[
N(j) = \sup_{k \in \Zd} |f(a_j^{-1}bk) - f(a_j^{-1}(x+bk))| .
\]
Since $|\det a_j^{-1}b| \# K(j)$ is bounded by the volume of
$F_2$, for all large $j$, and $N(j) \to 0$ by uniform continuity
of $f$, we deduce $R_j(x) \to 0$. This finishes the proof.

\subsection{\bf Proof of Theorem~\ref{surjectivity} --- synthesis
onto $L^p$} \ \label{surjectivity_proof}

Take $\lambda \in \C$ and $\sigma$ to be as in hypothesis
\eqref{bite}, and choose $\sigma^\prime \in (\sigma,1)$, so that
\[
\| \lambda P\psi(b \,  \cdot) - 1 \|_{L^p(\cube)}^p = \sigma <
\sigma^\prime .
\]
Take the analyzer to be $\phi = |b\cube|^{-1} \charfn_{b\cube}$, a
normalized indicator function which has constant periodization
$P|\phi| \equiv 1$.

Then for each $f \in L^p$,
\[
\| S_j (\lambda T_j f) - f \|_p^p \leq \sigma^\prime \| f \|_p^p
\]
for some $j > 0$, by Theorem~\ref{sample} applied to the function
$\lambda \psi$ (instead of to $\psi$). The coefficient sequence
here satisfies
\[
\| \lambda T_j f \|_{\ell^p(\Zd)} \leq |\lambda| |\det b|^{1-1/p}
\| f \|_p
\]
by formula \eqref{coeffexplicit} in the proof of
Theorem~\ref{continuityanal}.

Thus the open mapping theorem in Appendix~\ref{banachapp} says $S :
\ell^p(\Z_+ \times \Zd) \to L^p$ is open and surjective, and that
there exists $c \in \ell^p(\Z_+ \times \Zd)$ with $Sc = f$ and $\| c
\|_{\ell^p(\Z_+ \times \Zd)} \leq (1-\sigma^\prime)^{-1/p} |\lambda|
|\det b|^{1-1/p} \| f \|_p$. This completes the proof.

\subsection{\bf Examples and counterexamples for hypothesis \eqref{bite}} \ \label{tachev_proof}

Consider the function $\psi(x) = x^{-\beta} \charfn_{[0,1)}$ in one
dimension. We will show that if $\beta \in [2/(p+1),1/p)$ then $\| 1
- \lambda \psi \|_{L^p[0,1]} \geq 1$ for all $\lambda \in \C$, so
that hypothesis \eqref{bite} with $b=1$ fails for this function.
This counterexample for $p \in (0,1)$ is due to Tachev \cite{T95}.

Our proof below is different from Tachev's. It yields also the
positive result that the parameter range is sharp for $p \in (0,1]$:
if $\beta \in (0,2/(p+1))$ then $\| 1 - \lambda \psi \|_{L^p[0,1]} <
1$ for all small $\lambda>0$. Tachev stated this for $\beta=1$.

Assume $p \in (0,1)$ and $\beta \in [2/(p+1),1/p)$. To prove
Tachev's counterexample, we need only consider real, positive values
$\lambda>0$, since $\psi \geq 0$. Then after replacing $\lambda$ by
$t^{-\beta}$ and defining
\[
F(t) = \int_0^1 |1-(tx)^{-\beta}|^p \, dx ,
\]
we see we would like to prove $F(t) > 1$ for all $t>0$. A change of
variable gives
\begin{equation} \label{Feq1}
F(t) = \frac{1}{t} \int_0^t |1-y^{-\beta}|^p \, dy ,
\end{equation}
and hence $F$ is decreasing for $0 < t \leq 1$ because it equals the
mean value over the interval $[0,t]$ of the decreasing positive
function $(y^{-\beta}-1)^p$. For $t \geq 1$ we have
\begin{align}
F(t)
& = 1 + \frac{1}{t} \int_0^t ( |1-y^{-\beta}|^p - 1 ) \, dy \label{Feq} \\
& > 1 + \frac{1}{t} \int_0^\infty ( |1-y^{-\beta}|^p - 1 ) \, dy
\notag
\end{align}
because $|1-y^{-\beta}|^p - 1 < 0$ when $y > 1$. Note the last
integral converges near infinity because $\beta>1$. We now show this
integral is nonnegative. Indeed
\begin{align*}
& \int_0^\infty ( |1-y^{-\beta}|^p - 1 ) \, dy \\
& = \int_0^1 p(y^{-\beta} - 1)^{p-1} \beta y^{-\beta-1} y \, dy -
\int_1^\infty p(1 - y^{-\beta})^{p-1} \beta y^{-\beta-1} y \, dy \\
& \qquad \qquad \qquad \text{by parts on $(0,1)$ and $(1,\infty)$, using that $\beta p < 1$ and $\beta>1$,} \\
& = \beta p \int_0^1 (1 - y^\beta)^{p-1} y^{-\beta p} [1 - y^{\beta
(p+1)-2} ] \, dy
\end{align*}
by putting $y \mapsto y^{-1}$ in the second integral. Clearly the
integrand is nonnegative in this last integral, since $\beta (p+1)
\geq 2$, and this implies $F(t)>1$ for all $t>0$, as we wanted.

For the positive result when $\beta \in (1,2/(p+1))$, we simply use
\eqref{Feq} to prove
\[
F(t) = 1 + \frac{1}{t} \int_0^\infty ( |1-y^{-\beta}|^p - 1 ) \, dy
+ o \left( \frac{1}{t} \right) \qquad \text{as $t \to \infty$,}
\]
and then note the integral is negative by the calculations above.
Hence $F(t)<1$ for all large $t$, which shows $\| 1 - \lambda \psi
\|_{L^p[0,1]} < 1$ for all small $\lambda>0$. Next, when $\beta \in
(0,1]$ we observe that for each fixed $t>1$, the formula
\eqref{Feq1} for $F(t)$ is increasing with respect to $\beta$;
therefore $F(t)<1$ for all large $t>1$ by the case $\beta \in
(1,2/(p+1))$ just treated.

Lastly, the positive result when $p=1$ and $\beta \in (0,1)$ can
easily be proved directly.

\subsection{\bf Proof of Proposition~\ref{bitten} --- sufficient conditions} \ \label{bitten_proof}

Write $\eta(x) = P\psi(bx)$, so that $\eta \in L^p(\cube)$. Our
goal is to prove $\| 1 - \lambda \eta \|_{L^p(\cube)} < 1$ for
some $|\lambda| < 1$.

\emph{Part (a).} Suppose $\psi \in L^1$ and $\int_\Rd \psi \, dx
\neq 0$. Then $\eta \in L^1(\cube)$ with $\int_\cube \eta \, dx =
\int_\Rd \psi \, dx \neq 0$. We will use below the elementary
inequality that
\[
|1-z| \leq 1 - \Real z + A|z|^2 , \qquad |z| \leq 2^{-2/3},
\]
for some positive constant $A$. Given $0< |\lambda| \leq 1/2$, put
\[
E(\lambda) = \{ x \in \cube : |\eta(x)| \leq |\lambda|^{-1/3} \} ,
\]
and notice that on $E(\lambda)$ we have $|\lambda \eta| \leq
|\lambda|^{2/3} \leq 2^{-2/3}$. Then we see
\begin{align*}
\| 1 - \lambda \eta \|_{L^p(\cube)} & \leq \| 1 - \lambda \eta
\|_{L^1(\cube)} \qquad \text{by Jensen's inequality, since $p \leq 1$,} \\
& \leq \int_{E(\lambda)} (1 - \Real (\lambda \eta(x)) + A|\lambda
\eta(x)|^2) \, dx + \int_{\cube \setminus E(\lambda)} (1 + |\lambda \eta(x)|) \, dx \\
& = 1 - \Real \lambda \int_\cube \eta(x) \, dx + A|\lambda|^{4/3} +
o(|\lambda|) \qquad \text{as $|\lambda| \downarrow 0$,}
\end{align*}
since $E(\lambda) \uparrow \cube$. Thus we have only to choose
$\lambda$ with $\lambda \int_\cube \eta(x) \, dx > 0$ and
$|\lambda|$ sufficiently small, in order to obtain $\| 1 - \lambda
\eta \|_{L^p(\cube)} < 1$ as desired.

Notice the above proof works for $p=1$ as well.

\emph{Part (b).} Suppose $\psi \in L^1$ (so that $\eta \in
L^1(\cube)$) and that $\eta \not \equiv 0$ is real valued and
bounded above. (When $\eta$ is bounded below, just change
$\lambda$ to $-\lambda$ in what follows.) Suppose $\int_\Rd \psi
\, dx = 0$, so that $\int_\cube \eta \, dx = 0$. Then
\begin{align*}
\| 1 - \lambda \eta \|_{L^p(\cube)} & < \| 1 - \lambda \eta
\|_{L^1(\cube)} && \text{by Jensen's inequality} \\
& = \int_\cube (1 - \lambda \eta(x)) \, dx && \text{for all small $\lambda>0$, since $1-\lambda \eta>0$,} \\
& = 1 .
\end{align*}
Jensen's inequality is strict here because $p<1$ and $\lambda \eta$
is nonconstant (indeed, $\eta$ has mean value zero but is not
identically zero).

\emph{Part (c).} Assume $\eta \in L^2(\cube)$ satisfies $p
\int_\cube |\eta|^2 \, dx < (2-p) \left| \int_\cube \eta^2 \, dx
\right|$. Without loss of generality we can assume $\int_\cube
\eta^2 \, dx>0$, by multiplying $\eta$ with a suitable complex
constant. Then our assumption is equivalent to $\int_\cube (\Imag
\eta)^2 \, dx < (1-p) \int_\cube (\Real \eta)^2 \, dx$, so that
\begin{equation} \label{realim}
\alpha^{-1} \int_\cube (\Imag \eta)^2 \, dx < \alpha (1-p)
\int_\cube (\Real \eta)^2 \, dx ,
\end{equation}
for some $\alpha<1$ sufficiently close to $1$.

We will use below the binomial approximation that
\[
|1-z|^p \leq 1 - p \Real z + \alpha \frac{p(p-1)}{2} (\Real z)^2 +
\alpha^{-1} \frac{p}{2} (\Imag z)^2 , \qquad |z| \leq B,
\]
where the small positive constant $B$ depends on $p \in (0,1)$ and
$\alpha \in (0,1)$. Putting
\[
F(\lambda) = \{ x \in \cube : |\eta(x)| \leq B/|\lambda| \} ,
\]
we deduce that
\begin{align*}
& \| 1 - \lambda \eta \|_{L^p(\cube)}^p \\
& \leq \int_{F(\lambda)}
\left( 1 - p \Real (\lambda \eta(x)) + \alpha \frac{p(p-1)}{2}
(\Real \lambda \eta(x))^2 + \alpha^{-1} \frac{p}{2} (\Imag \lambda
\eta(x))^2 \right) dx \\
& \qquad \qquad + \int_{\cube \setminus F(\lambda)} (1 + |\lambda
\eta(x)|^p) \, dx .
\end{align*}
Averaging over $\lambda$ and $-\lambda$ (and noting
$F(-\lambda)=F(\lambda)$) gives for small $\lambda>0$ that
\begin{align*}
& \frac{1}{2} \left( \| 1 - \lambda \eta \|_{L^p(\cube)}^p + \|
1 + \lambda \eta \|_{L^p(\cube)}^p \right) \\
& \leq \int_{F(\lambda)} \left( 1 + \alpha \frac{p(p-1)}{2} (\Real
\lambda \eta(x))^2 + \alpha^{-1} \frac{p}{2} (\Imag \lambda
\eta(x))^2 \right) dx \\
& \qquad \qquad \qquad \qquad + \int_{\cube \setminus F(\lambda)} (1
+ |\lambda \eta(x)|^p) \, dx \\
& = 1 + \lambda^2 \int_\cube \left( \alpha \frac{p(p-1)}{2} (\Real
\eta(x))^2 + \alpha^{-1} \frac{p}{2} (\Imag \eta(x))^2 \right) \, dx
+ o(\lambda^2) \qquad \text{as $\lambda \downarrow 0$,}
\end{align*}
where in the final step we used that $F(\lambda) \uparrow \cube$ and
that $|\lambda \eta|^p < B^{p-2} |\lambda \eta|^2$ on $\cube
\setminus F(\lambda)$ (because $|\lambda \eta|/B > 1$ there).

Hence from \eqref{realim} we conclude
\[
\frac{1}{2} \left( \| 1 - \lambda \eta \|_{L^p(\cube)}^p + \| 1 +
\lambda \eta \|_{L^p(\cube)}^p \right) < 1
\]
for all small $\lambda>0$, and thus by choosing either $\lambda$ or
$-\lambda$ we complete the proof.

\subsubsection*{Aside.} Our proofs of parts (a) and (c)
are modifications of Filippov and Oswald \cite[Lemma~1]{FO95}. They
treated only real-valued functions $\eta$, which in the proof of
part (c) above means they could choose $B=1/2$ and $\alpha$
sufficiently close to $0$, whereas we must choose $\alpha$
sufficiently close to $1$ and then take $B$ sufficiently close to
$0$.

\subsection{\bf Proof of Theorem~\ref{undersynth} ---
more synthesis onto $L^p$} \ \label{undersynth_proof}

Let $\beta \in \N$ and consider the periodization of $\psi$ with
respect to the translation matrix $b\beta$, that is,
\[
P_{b\beta}\psi(x) = |\det b\beta| \sum_{k \in \Zd} \psi(x-b\beta k)
.
\]
After rescaling, we see $P_{b\beta}\psi(b\beta x)$ is integrable on
the cube $\cube$ (since $\psi \in L^1$), and has Fourier
coefficients $\widehat{\psi}(m(b\beta)^{-1})$ for $m \in \Zd$ (row
vectors).

We claim $P_{b\beta} \psi$ is nontrivial for some $\beta$. For if
$P_{b\beta} \psi =  0$ a.e. for each $\beta$, then the Fourier
coefficients are zero too, so that
$\widehat{\psi}(m(b\beta)^{-1})=0$ for all $m \in \Zd$ and all
$\beta \in \N$. Then continuity of $\widehat{\psi}$ forces
$\widehat{\psi} \equiv 0$, contradicting the hypothesis that $\psi
\not \equiv 0$.

So fix a $\beta$ value for which $P_{b\beta} \psi$ is nontrivial.
This periodization is real-valued (since $\psi$ is real valued), and
is bounded below since $\psi_-$ is bounded and has compact support.
Thus Proposition~\ref{bitten}(a) or \ref{bitten}(b) applies, and
says $\| \lambda P_{b\beta} \psi(b\beta \cdot) - 1 \|_{L^p(\cube)} <
1$ for some $\lambda$.

Theorem~\ref{surjectivity} then provides a constant $C$ such that
for each $f \in L^p$ there is a sequence $\tilde{c} \in \ell^p(\Z_+
\times \Zd)$ with $\| \tilde{c} \|_{\ell^p(\Z_+ \times \Zd)} \leq C
\| f \|_p$ and
\[
f = \sum_{j>0} \sum_{k \in \Zd} \tilde{c}_{j,k} |\det a_j|^{1/p}
\psi(a_j x - b\beta k) .
\]
That is, $f=Sc$ where $c_{j,k}=\tilde{c}_{j,\beta^{-1}k}$ if $k \in
\beta \Zd$ and $c_{j,k}=0$ otherwise. The theorem follows, since $c$
and $\tilde{c}$ have the same $\ell^p$-norm.

\subsection{\bf Proof of Corollary~\ref{atomicdomain} ---
affine atomic decomposition of $L^p(\Omega)$} \
\label{atomicdomain_proof}

The ``$\leq$'' direction of the Corollary follows straight from
Proposition~\ref{continuitysynth}.

For the ``$\geq$'' direction, first define
\[
{\mathcal L} = \{ c \in \ell^p(\Z_+ \times \Zd) : \text{$c$ is
adapted to $\Omega$ and $\psi$} \} .
\]
Clearly ${\mathcal L}$ is a closed subspace of $\ell^p(\Z_+ \times
\Zd)$, and hence is a complete metric space under the
$\ell^p$-metric. Take $\phi = |b\cube|^{-1} \charfn_{b\cube}$.

Consider an $f \in L^p$ that is continuous and compactly supported
in $\Omega$. We claim the sequence $T_j f$ belongs to ${\mathcal
L}$, for each large $j$, or more precisely, that the sequence $c$
equalling $T_j f$ at level $j$ and zero at all other levels belongs
to ${\mathcal L}$. To see this, just notice
\[
\text{$\spt(\psi_{j,k}) \subset \Omega$ whenever $\langle \Theta f ,
\phi(a_j \, \cdot - bk) \rangle \neq 0$ and $k \in \Zd$,}
\]
for all large $j$, because the support of $f$ lies at some positive
distance from the boundary of $\Omega$, and $\psi$ and $\phi$ have
compact support and $\| a_j^{-1} \| \to 0$. Thus $T_j f$ belongs to
${\mathcal L}$.

The proof of Theorem~\ref{surjectivity} now applies word-for-word,
except with $\ell^p(\Z_+ \times \Zd)$ replaced by ${\mathcal L}$.
Admittedly we have verified the hypotheses of the open mapping
theorem only for the dense class of continuous $f$ having compact
support in $\Omega$, but a dense class in $L^p(\Omega)$ is enough,
by the remark in Appendix~\ref{banachapp}.

The conclusion of Theorem~\ref{surjectivity} with $c \in {\mathcal
L}$ gives the ``$\geq$'' direction of Corollary~\ref{atomicdomain}.

\subsection{\bf Proof of Corollary~\ref{coeffanal} --- analysis metric for $L^p$} \
\label{coeffanal_proof}

By formula \eqref{coeffexplicit} in the proof of
Theorem~\ref{continuityanal},
\[
\sup_j \| T_j f \|_{\ell^p} \leq |\det b|^{1-1/p}
\|P|\phi|\|_\infty^{1/p} \|f\|_p.
\]
To prove the other direction of the metric equivalence, choose a
synthesizer $\psi = |b\cube|^{-1} \charfn_{b\cube}$ that has
constant periodization $P\psi \equiv 1$. Then by
Theorem~\ref{sample}, $S_j T_j f \to f$ in $L^p$ as $j \to \infty$.
Therefore Proposition~\ref{continuitysynth} (bounded synthesis)
implies that
\[
\| f \|_p \leq \sup_j \| S_j T_j f \|_p \leq \| \psi \|_p \sup_j
\| T_j f \|_{\ell^p(\Zd)} .
\]

\subsection{\bf Proof of Theorem~\ref{inject} ---
$p-$Riesz basis at scale $j$} \ \label{inject_proof}

By a simple rescaling, we can suppose $a_j=I$ is the identity
matrix.

To prove injectivity, take $s \in \ell^p(\Zd)$ and suppose $S_j s =
0$, or
\[
\sum_{k \in \Zd} s_k \psi(x-bk) = 0 \quad \text{a.e.}
\]
We will show $s=0$.

Note the series $\sum_{k \in \Zd} s_k \psi(x-bk)$ converges
absolutely a.e., because $s \in \ell^p \subset \ell^\infty$ and
\[
(\sum_{k \in \Zd} |\psi(x-bk)|)^p \leq \sum_{k \in \Zd}
|\psi(x-bk)|^p \in L^1(b\cube) .
\]
(Here we use that $p \in (0,1]$.) Let $\xi \in \Rd$ and multiply
the series by $e^{2\pi i \xi x}$ (where $\xi \in \Rd$ is
arbitrary) to obtain
\[
\sum_{k \in \Zd} s_k e^{2\pi i \xi bk} e^{2\pi i \xi (x-bk)}
\psi(x-bk) = 0 \quad \text{a.e.}
\]
Replace $x$ by $x-b\ell$ and sum over $\ell \in \Zd$ to obtain that
\[
\sum_{k \in \Zd} s_k e^{2\pi i \xi bk} |\det b| \sum_{\ell \in \Zd}
e^{2\pi i \xi (x-b\ell -bk)} \psi(x-b\ell -bk) = 0 \quad
\text{a.e.,}
\]
with the double series converging absolutely a.e.\ because $s \in
\ell^p \subset \ell^1$. Thus
\[
\sum_{k \in \Zd} s_k e^{2\pi i \xi bk} \cdot P(e_\xi \psi)(x) = 0
\quad \text{a.e.}
\]
By hypothesis there is a set of positive measure on which $P(e_\xi
\psi)(x) \neq 0$, and hence $\sum_{k \in \Zd} s_k e^{2\pi i \xi
bk} = 0$ for each $\xi \in \Rd$. Since $s \in \ell^p \subset
\ell^1$ we conclude $s_k=0$ for all $k$, or $s=0$, so that $S_j$
is injective.

Now suppose in addition that $\psi$ has compact support. We will
prove the $p$-Riesz basis condition by following almost exactly
the work of R.-Q. Jia \cite[\S3]{J98}. Our proof does present one
new idea: whereas Jia restricted his $\psi \in L^p$ to be a
distribution, so that he could work with its Fourier transform, we
avoid any such restriction by working directly with the
periodization hypothesis.

Define
\[
\psi^{(\ell)}(x) =
\begin{cases}
\psi(x+b\ell) , & x \in b\cube , \\
0 , & \text{otherwise,}
\end{cases}
\]
so that $\psi^{(\ell)}$ gives the value of $\psi$ on
$b(\ell+\cube)$, translated to $b\cube$. Obviously $\psi$ can be
reconstructed by summing up the pieces:
\begin{equation} \label{reconstruct}
\psi = \sum_{\ell \in \Zd} \psi^{(\ell)}(\cdot - b\ell) .
\end{equation}
Only finitely many of the $\psi^{(\ell)}$ are nontrivial, since
$\psi$ has compact support, and so we can choose a maximal
collection of them that are linearly independent in $L^p(b\cube)$.
Denote this collection by $\{ \psi^{(m)} : m \in M \}$ for some
finite index set $M \subset \Zd$.

For later use, write $v=\{ v_m \}_{m \in M}$ for an arbitrary
complex sequence supported on $M$, and observe that the function
$f(v) = \| \sum_{m \in M} v_m \psi^{(m)} \|_{L^p(b\cube)}$ is
continuous on the unit $p$-sphere $\{ v : \| v \|_{\ell^p(M)} = 1
\}$. Clearly $f$ cannot equal zero anywhere on this sphere,
because the $\psi^{(m)}$ are linearly independent. Hence $f$
attains a positive minimum value $C=C(\psi,b,p)$ on the unit
$p$-sphere. (Finiteness of the index set $M$ is used here to
ensure compactness of the unit sphere, and hence existence of a
minimum for $f$.) Thus
\begin{equation} \label{homogeq}
\| \sum_{m \in M} v_m \psi^{(m)} \|_{L^p(b\cube)} \geq C (\sum_{m
\in M} |v_m|^p)^{1/p} , \qquad v_m \in \C,
\end{equation}
by homogeneity.

Each $\psi^{(\ell)}$ can be expressed as a linear combination
\begin{equation} \label{lincombeq}
\psi^{(\ell)} = \sum_{m \in M} t_{\ell,m} \psi^{(m)} , \qquad \ell
\in \Zd
\end{equation}
for some coefficients $t_{\ell,m}$. Substituting this into the
reconstruction formula \eqref{reconstruct} gives
\[
\psi(x) = \sum_{\ell \in \Zd} \sum_{m \in M} t_{\ell,m}
\psi^{(m)}(x-b\ell) .
\]
Hence
\begin{align*}
S_j s(x)
& = \sum_{k \in \Zd} s_k \psi(x-bk) \qquad \text{(recalling that $a_j=I$)} \\
& = \sum_{k \in \Zd} \sum_{\ell \in \Zd} \sum_{m \in M} s_k
t_{\ell,m} \psi^{(m)}(x-bk-b\ell) \\
& = \sum_{\ell \in \Zd} \sum_{m \in M} (s*t_m)_\ell \,
\psi^{(m)}(x-b\ell)
\end{align*}
by shifting the index $\ell \mapsto \ell-k$ and defining a
sequence $t_m = \{ t_{\ell,m} \}_{\ell \in \Zd}$, for each $m \in
M$. Convergence of the above multiple series is clear, because
each sequence $t_m$ has only finitely many nonzero entries
$t_{\ell,m}$ (noting $\psi^{(\ell)}$ is identically zero for all
large $|\ell|$).

Since $\psi^{(m)}$ equals zero outside the cube $b\cube$, we
deduce
\[
S_j s(x) = \sum_{m \in M} (s*t_m)_\ell \, \psi^{(m)}(x-b\ell) ,
\qquad x \in b(\ell + \cube), \quad \ell \in \Zd .
\]
Therefore
\begin{align}
\| S_j s \|_p^p
& = \sum_{\ell \in \Zd} \| S_j s \|_{L^p(b(\ell+\cube))}^p \notag \\
& = \sum_{\ell \in \Zd} \| \sum_{m \in M} (s*t_m)_\ell
\, \psi^{(m)} \|_{L^p(b\cube)}^p \notag \\
& \geq C \sum_{\ell \in \Zd} \sum_{m \in M} |(s*t_m)_\ell|^p
\qquad \text{by \eqref{homogeq}} \notag \\
& = C \sum_{m \in M} \| s*t_m \|_{\ell^p(\Zd)}^p .
\label{lowerest}
\end{align}

We must still bound the norm of $s*t_m$ from below in terms of the
norm of $s$. To help achieve this, consider the trigonometric
polynomial $\tau_m(\xi)=\sum_{\ell \in \Zd} t_{\ell,m} e^{2\pi
i\xi \ell}$. For each $\xi \in \Rd$, our periodization hypothesis
guarantees that
\begin{align*}
0 \not \equiv P(e_\xi \psi)(x)
& = |\det b| \sum_{\ell \in \Zd} e^{2\pi i \xi(x+b\ell)} \psi(x+b\ell) \\
& = |\det b| e^{2\pi i \xi x} \sum_{\ell \in \Zd} e^{2\pi i \xi b\ell} \psi^{(\ell)}(x) \qquad \text{for $x \in b\cube$} \\
& = |\det b| e^{2\pi i \xi x} \sum_{m \in M} \tau_m(\xi b)
\psi^{(m)}(x)
\end{align*}
by substituting \eqref{lincombeq}. We deduce that at least one of
the values $\tau_m(\xi b), m \in M$, must be nonzero. Hence
$\sum_{m \in M} |\tau_m(\xi)|^2 > 0$ for all $\xi$, and so the
reciprocal function
\[
\upsilon(\xi)= ( \sum_{m \in M} |\tau_m(\xi)|^2 )^{-1} .
\]
is well defined, smooth and $\Zd$-periodic. Write $u_\ell$ for its
Fourier coefficients: $\sum_{\ell \in \Zd} u_\ell e^{2\pi i \xi
\ell} = \upsilon(\xi)$. These Fourier coefficients decay rapidly,
since $\upsilon$ is smooth. And writing
$\widetilde{t}_{\ell,m}=\overline{t_{-\ell,m}}$, we have from
\eqref{lowerest} the estimate
\begin{align*}
\| S_j s \|_p^p & \geq C \frac{\sum_{m \in M} \| s*t_m
\|_{\ell^p(\Zd)}^p \| \widetilde{t}_m * u
\|_{\ell^p(\Zd)}^p}{\max_{m \in
M} \| \widetilde{t}_m * u \|_{\ell^p(\Zd)}^p} \\
& \geq C \| \sum_{m \in M} s*t_m * \widetilde{t}_m
* u \|_{\ell^p(\Zd)}^p \\
& = C \| s \|_{\ell^p(\Zd)}^p
\end{align*}
since $\sum_{m \in M} t_m * \widetilde{t}_m * u = \delta$, as one
can check by taking the Fourier series: $\sum_{m \in M} \tau_m
\overline{\tau_m} \upsilon = 1$.

Thus we have proved the lower Riesz estimate for the theorem. The
upper estimate $\| S_j s \|_p \leq C \| s \|_{\ell^p(\Zd)}$ is
immediate from Proposition~\ref{continuitysynth}. Now the range
$S_j(\ell^p(\Zd))$ must be complete in $L^p$, as one sees by
considering Cauchy sequences in the range and using \eqref{inject1},
and so the range is closed.

\subsection{\bf Proof of Theorem~\ref{properties} ---
most $L^p$ functions are surjective affine synthesizers} \
\label{properties_proof}

[Density.] The class ${\mathcal S}^p$ contains every bounded
function $\psi$ with compact support and nonzero integral, because
every such function satisfies the hypotheses of
Proposition~\ref{bitten}(a) and hence of Theorem~\ref{surjectivity}.
These bounded functions are dense in $L^p$, and hence ${\mathcal
S}^p$ is dense in $L^p$.

[Openness.] Take $\psi \in {\mathcal S}^p$. Then $S_\psi$ is a
continuous linear mapping of the $F$-space $\ell^p(\Z_+ \times \Zd)$
onto the $F$-space $L^p$, so that $S_\psi$ is open by
\cite[Corollary 2.12]{R91}. Hence $A>0$ exists such that for each $f
\in L^p$, a sequence $c \in \ell^p(\Z_+ \times \Zd)$ exists
satisfying $S_\psi c=f$ and
\[
d_{\ell^p}(0,c) \leq A d_p(0,f) .
\]
We claim ${\mathcal S}^p$ contains the $L^p$-ball of radius $1/A$
centered at $\psi$, from which it follows that ${\mathcal S}^p$ is
open in $L^p$.

So suppose $\psi_1 \in L^p$ with $d_p(\psi_1,\psi) = \delta/A$ for
some $\delta \in (0,1)$. Then
\begin{align*}
d_p(S_{\psi_1}c,f)
& = \| S_{(\psi_1-\psi)}c \|_p^p \\
& \leq \| \psi_1-\psi \|_p^p \| c \|_{\ell^p(Z_+ \times
\Zd)}^p && \text{by Proposition~\ref{continuitysynth}} \\
& \leq \frac{\delta}{A} A d_p(0,f) = \delta d_p(0,f) && \text{by
construction above.}
\end{align*}
Now the open mapping theorem in Appendix~\ref{banachapp} guarantees
that $S_{\psi_1}$ maps onto $L^p$, so that $\psi_1 \in {\mathcal
S}^p$ as desired.

[Path connectedness.] First we show path connectedness of
\begin{equation} \label{subclass}
\{ \psi \in L^p : \| P\psi(b \, \cdot) - 1 \|_{L^p(\cube)} < 1 \} ,
\end{equation}
which is a subset of ${\mathcal S}^p$ by Theorem~\ref{surjectivity}.
Consider the linear path
\[
\psi_t = (1-t) \psi + t |b\cube|^{-1} \charfn_{b\cube} , \qquad t
\in [0,1] ,
\]
which connects $\psi$ to the normalized indicator function
$|b\cube|^{-1} \charfn_{b\cube}$. This normalized indicator function
has periodization identically equal to $1$, and so
\begin{align*}
\| P\psi_t(b \, \cdot) - 1 \|_{L^p(\cube)}
& = \| P\psi(b \, \cdot) - 1 \|_{L^p(\cube)} (1-t) \\
& \leq \| P\psi(b \, \cdot) - 1 \|_{L^p(\cube)} < 1 .
\end{align*}
This proves path connectedness of the collection \eqref{subclass},
as we wanted.

It follows immediately that the collection
\begin{equation} \label{bitelater}
\{ \psi \in L^p : \| \lambda P\psi(b \, \cdot) - 1 \|_{L^p(\cube)} <
1 \text{\ for some $\lambda \in \C , \lambda \neq 0$} \}
\end{equation}
is also path connected and lies in ${\mathcal S}^p$, because
$\lambda \psi$ belongs to the collection \eqref{subclass} and $\psi$
is path connected to $\lambda \psi$ within the collection
\eqref{bitelater}, through an obvious path of rescalings.

Now consider an arbitrary $\widetilde{\psi} \in {\mathcal S}^p$. By
openness, there exists an $L^p$-ball around $\widetilde{\psi}$ that
lies in ${\mathcal S}^p$. This ball contains some bounded function
$\psi$ having compact support and nonzero integral, and this $\psi$
belongs to the collection \eqref{bitelater} by
Proposition~\ref{bitten}(a). We can connect $\widetilde{\psi}$ to
$\psi$ by a path lying in the ball, and so path connectedness of
${\mathcal S}^p$ follows from path connectness of collection
\eqref{bitelater}.

\section*{\bf Acknowledgments}
John Benedetto sparked this research during the International
Conference on Harmonic Analysis and Applications, Villa de Merlo,
Argentina (2006), by asking me about $L^p$-affine synthesis for
$0<p<1$.

Qui Bui is my co-author on the papers \cite{bl1}--\cite{bl5}. The
many fruitful discussions we enjoyed while writing those works have
influenced the current paper as well.

\appendix

\section{\bf The open mapping theorem}
\label{banachapp}

The open mapping theorem for metric spaces was used in the following
form, in the proof of Theorem~\ref{surjectivity} (surjectivity of
the synthesis operator).
\begin{proposition} \label{banach}
Let $X$ and $Y$ be complete metric vector spaces with
translation-invariant metrics $d_X$ and $d_Y$ respectively. Suppose
$S : X \to Y$ is continuous and linear, take $\delta \in (0,1)$ and
$A>0$, and assume for each $y \in Y$ that some $x \in X$ exists with
\begin{equation} \label{iterate}
d_Y(Sx,y) \leq \delta d_Y(0,y) , \qquad d_X(0,x) \leq A d_Y(0,y) .
\end{equation}

Then $S$ is an open mapping, and $S(X)=Y$. Indeed, given $y \in Y$
there exists $x \in X$ with $Sx=y$ and $d_X(0,x) \leq
(1-\delta)^{-1} A d_Y(0,y)$.
\end{proposition}

\noindent \emph{Remark.} The hypothesis in Proposition~\ref{banach}
can be weakened to assume only for some \emph{dense} subset of
$y$-values that $x$ exists satisfying \eqref{iterate}, provided we
are prepared to replace $\delta$ in the conclusion of the
Proposition by $\delta^* \in (\delta,1)$ and $A$ by $A^* > A$.
\begin{proof}[Proof of Proposition~\ref{banach}]
Let $y_0 \in Y$. Choose $x_0 \in X$ according to \eqref{iterate}
with $y=y_0$. Let $y_1=y_0 - Sx_0$ and choose $x_1$ according to
\eqref{iterate} with $y=y_1$. Let $y_2=y_1 - Sx_1$, and continue
this process, obtaining $x_0,x_1,x_2,\ldots \in X$ and
$y_0,y_1,y_2,\ldots \in Y$ that satisfy
\begin{align}
y_{m+1} & = y_m - Sx_m , \label{banacheq1} \\
d_Y(0,y_{m+1}) & \leq \delta d_Y(0,y_m) , \label{banacheq2} \\
d_X(0,x_m) & \leq A d_Y(0,y_m) , \label{banacheq3}
\end{align}
for $m=0,1,2,\ldots$. (The lefthand side of \eqref{banacheq1} uses
the translation invariance of the $Y$-metric.)

Now define $x = \sum_{m=0}^\infty x_m$, which converges in the
complete, translation-invariant space $X$ because
\begin{align}
\sum_{m=0}^\infty d_X(0,x_m) & \leq A
\sum_{m=0}^\infty d_Y(0,y_m) && \text{by \eqref{banacheq3}} \notag \\
& \leq A d_Y(0,y_0) \sum_{m=0}^\infty \delta^m && \text{by \eqref{banacheq2}} \notag \\
& = \frac{A}{1-\delta} d_Y(0,y_0) . \label{banacheq4}
\end{align}
The continuity and linearity of $S$ imply that
\[
Sx = \sum_{m=0}^\infty Sx_m = \sum_{m=0}^\infty (y_m-y_{m+1}) = y_0
,
\]
by \eqref{banacheq1} and telescoping, since $y_{m+1} \to 0$ by
\eqref{banacheq2}. Because $y_0$ was arbitary, we have shown
$S(X)=Y$. Further, \eqref{banacheq4} shows
\[
d_X(0,x) \leq \frac{A}{1-\delta} d_Y(0,y_0) .
\]
It follows for all $r>0$ that $S(B_X(r)) \supset
B_Y((1-\delta)A^{-1}r)$, where ``$B$'' denotes an open ball, and
thus $S$ is an open mapping.
\end{proof}
\end{document}